\documentclass[leqno,12pt]{article}%
\usepackage{amsmath,amssymb}
\usepackage{dsfont}
\usepackage{amsfonts}%
\usepackage{color}

\newtheorem{Theorem}{\hspace{\parindent}\bf Theorem}[section]
\newtheorem{Lemma}{\hspace{\parindent}\bf Lemma}[section]
\newtheorem{Proposition}{\hspace{\parindent}\bf Proposition}[section]
\newtheorem{Corollary}{\hspace{\parindent}\bf Corollary}[section]

\newtheorem{Remark}{\hspace{\parindent}\bf Remark}[section]

\newcommand{\qed}{\hfill$\square$\vspace{0.3cm}}
\begin{document}

\title{\textbf{Regularity theory for singular nonlocal diffusion equations}}
\author{by\\Arturo de Pablo, Fernando Quir\'{o}s, and Ana Rodr\'{\i}guez}

\maketitle

\

\begin{abstract}
We prove continuity for bounded  weak solutions of a nonlinear nonlocal parabolic type equation associated to a Dirichlet form with a rough kernel. The equation is allowed to be singular  at the level zero, and solutions may  change sign.   If  the nonlinearity in the equation does not oscillate too much at the origin,  the solution is proved to be moreover H\"older continuous.

The results are new  even when the Dirichlet form is the one corresponding to the fractional Laplacian.
\end{abstract}


\vskip 8cm

\noindent{\makebox[1in]\hrulefill}\newline2010 \textit{Mathematics Subject
Classification.} 35R11, 
35B65, 
35K55. 
\newline\textit{Keywords and phrases.} Singular nonlocal diffusion, integral operators, De Giorgi methods.
regularity.

\newpage

\section{Introduction and main results}

\label{sect-introduction} \setcounter{equation}{0}

The aim of this paper is to prove regularity of bounded weak solutions $u$ to
\begin{equation}
\label{eq:main}
\partial_t \beta(u)+\mathcal{L} u=0,
\end{equation}
where $\mathcal{L}$ is a nonlocal operator associated to the bilinear Dirichlet form
\begin{equation*}
\label{quadratic-form}\mathcal{E}_J(u,v)=\frac12\int_{\mathbb{R}^{N}}%
\int_{\mathbb{R}^{N}}(u(x)-u(y))(v(x)-v(y))J(x,y)\,dxdy,
\end{equation*}
with a measurable kernel $J$ satisfying
\begin{equation}
\label{eq:kernel}
\tag{$\text{\rm H}_{J}$}
\left\{
\begin{array}{l}
J(x,y)\ge 0, \quad J(x,y)=J(y,x),\\[10pt]
\displaystyle
\frac{\mathds{1}_{\{|x-y|\le
1\}}}{\nu |x-y|^{N+\sigma}}\le J(x,y)\le\frac{\nu}{|x-y|^{N+\sigma}},
\end{array}
\right.
\qquad  \text{for a.e. }(x,y)\in
\mathbb{R}^{N}\times\mathbb{R}^N,
\end{equation}
for some constants $\sigma\in(0,2)$ and $\nu\ge1$. The bilinear form
is defined in $\mathcal{H}_J$, which  is the space of measurable functions with  $\mathcal{E}_J(u,u)<\infty$. Since we are not asking $J$ to be smooth outside the diagonal $x=y$, it is referred to in the literature as a \emph{rough kernel}.  For the smooth (outside the diagonal) kernel $J(x,y)=|x-y|^{-N-\sigma}$, the operator is a multiple of the well-known fractional Laplacian $(-\Delta)^{\sigma/2}$.

%
%

To be precise, $\mathcal{L}: \mathcal{H}_J\to \mathcal{H}_J'$ is the linear operator defined by $\langle \mathcal{L}u,\zeta\rangle =\mathcal{E}_J(u,\zeta)$ for any $\zeta\in C^\infty_{\rm c}(\mathbb{R}^N)$. Thus, $u$ is a weak solution to equation \eqref{eq:main} if
\begin{equation}
\label{eq:main.weak}
\left\{\begin{array}{l}
u\in L^2_{\rm loc}\big(\mathbb{R}_+:\mathcal{H}_J\big), \quad \beta(u)\in C(\mathbb{R}_+: L^{1}(\mathbb{R}^{N})),
\\[10pt]
\displaystyle\int_0^\infty\int_{\mathbb{R}^{N}}\beta(u)\partial_{t} \zeta\,dxdt-\int_0^\infty\mathcal{E}_J(u,\zeta
)\,dt=0\quad \text{for all }\zeta\in C^\infty_c(\mathbb{R}^N\times\mathbb{R}_+).
\end{array}
\right.
\end{equation}
Under assumptions \eqref{eq:kernel}, existence and uniqueness of a weak solution to \eqref{eq:main} that is moreover bounded, with a prescribed initial value $\beta(u(\cdot,0))\in L^1(\mathbb{R}^N)\cap L^\infty(\mathbb{R}^N)$, are proved in~\cite{PQR} whenever $\beta$ is continuous and nondecreasing.

In this paper the nonlinearity $\beta$  is assumed to satisfy moreover
\begin{equation}\label{eq:cond.beta0}\tag{${\rm H}_0$}
\beta\in C^1(\mathbb{R}),\qquad\beta(0)=0,\qquad\beta'(s)>0\quad\text{for } s\ne0.
\end{equation}
Notice that  we are allowing $\beta'(0)$ to be zero. If this is the case,
the equation is \emph{singular}, since the \lq\lq diffusion'' coefficient $1/\beta'(u)$ becomes singular at the level $u=0$.

Since we are dealing with bounded solutions, when proving regularity for a given solution $u$ we may replace $\beta$ linearly for $|s|>\|u\|_\infty$, so we always assume $\beta'$ bounded.

We prove that bounded weak solutions to equation~\eqref{eq:main} are continuous if the equation is \emph{not too singular}, namely, if $\beta'$ is bounded from below by some power,
\begin{equation}\label{eq:cond.beta1}\tag{${\rm H}_1$}
\beta'(s)\le\mathcal{M}_0\;\mbox{ for } s\in\mathbb{R},\quad\beta'(s)\ge \mathfrak{m}_0|s|^{p-1}\;\mbox{ for } |s|\le \|u\|_\infty,\quad \mathfrak{m}_0>0,\,p>1.
\end{equation}

\begin{Theorem}
\label{th:regularity} Let $J$ and $\beta$ satisfy respectively~\eqref{eq:kernel} and~\eqref{eq:cond.beta0}--\eqref{eq:cond.beta1}. Then bounded weak solutions to equation~\eqref{eq:main} are continuous in $\mathbb{R}^{N}\times\mathbb{R}_+$.
\end{Theorem}
We remark that our solutions may change sign, which introduces an extra main difficulty precisely when $\beta'(0)=0$.

To get further regularity, we will require in addition that $\beta'$  does not oscillate too much at the origin. More precisely, we will assume that  there exists a slowly varying function $h:(0,\|u\|_\infty)\to(0,\infty)$,  such that
\begin{equation}\label{eq:cond.beta2}\tag{${\rm H}_2$}
\mathfrak{m}_1\le \dfrac{|s|^{1-q}\beta'(s)}{h(|s|)}\le \mathcal{M}_1,\qquad 0<|s|\le \|u\|_\infty,\quad \mathfrak{m}_1,\,\mathcal{M}_1> 0,\,q\ge1.
\end{equation}
Roughly speaking this means that $\beta'$ has the order of a power at the origin, possibly perturbed by a lower order (bounded or unbounded) coefficient.
We recall that, according to Karamata~\cite{Karamata}, a measurable function $h: (0,1)\to (0,\infty)$ varies slowly (at zero) if
$$
\lim_{\tau\to0^+}\frac{h(\lambda \tau)}{h(\tau)}=1\quad \text{ for every } \lambda>0.
$$

As examples of slowly varying functions we have $h(\tau)=c$, $h(\tau)=c(\log(1/\tau))^d$ or $h(\tau)=c(\log(\log(e/\tau)))^d$, $d\in\mathbb{R}$.
See a complete account of the theory in~\cite{Bingham-Goldie-Teugels}. Observe that this allows us to include  the case of weak singular (almost linear) diffusion $\beta(s)=s(\log(2/|s|))^{-1}$ for $|s|\le1$.

\begin{Theorem}
\label{th:holder.continuity} If in addition to the conditions of Theorem~\ref{th:regularity} we also have \eqref{eq:cond.beta2}, then  bounded weak solutions to equation~\eqref{eq:main} are H\"older continuous  at every point.
\end{Theorem}

It turns out that the H\"older constants and exponents are uniform in sets where $u\ge\delta>0$. As a corollary we get a better result when the solution is positive.
\begin{Corollary}
\label{th:regularity>0} Under the conditions of Theorem~\ref{th:holder.continuity}, positive solution are uniformly H\"older continuous in every compact set of $\mathbb{R}^N\times\mathbb{R}_+$.
\end{Corollary}
We recall that if  $u(\cdot,0)\ge0$ in $\mathbb{R}^N$ then $u(\cdot,t)>0$ in $\mathbb{R}^N$ for every $t>0$; cf~\cite{VPQR}.

The assumption \lq\lq$u$ bounded'' in the previous regularity results is not a big restriction. In fact  solutions become immediately bounded provided some integrability condition is imposed on the initial value, see \cite{PQRV2,VPQR} for the fractional Laplacian case.

Condition~\eqref{eq:cond.beta0} guarantees the existence of an inverse $\varphi=\beta^{-1}$ for $\beta$. The function $v=\beta(u)$, which inherits the regularity that we have obtained for $u$, satisfies the \emph{nonlocal filtration} equation
\begin{equation}
\label{eq:elefi}
\partial_t v+\mathcal{L}\varphi(v)=0
\end{equation}
in a weak sense. Thus, when $\beta(u)=|u|^{\frac1m-1}u$, with $m\in(0,1)$, we get regularity for bounded weak solutions to the nonlocal fast diffusion type equation
\begin{equation*}
\label{eq:FFDE}
\partial_tv+\mathcal{L}(|v|^{m-1}v)=0,\qquad 0<m<1,
\end{equation*}
which was not known even for the fractional Laplacian case. Moreover, the regularity result obtained in Corollary \ref{th:regularity>0} is enough to prove that bounded weak solutions  of $\partial_t v+(-\Delta)^{\sigma/2}\varphi(v)=0$ with a sign are classical solutions provided $\varphi$ is regular enough; see~\cite{VPQR}.

\textsc{Precedents. } When $\mathcal{L}=(-\Delta)^{\sigma/2}$ the problem can be transformed into a local one by means of the extension technique introduced by Caffarelli and Silvestre in~\cite{Caffarelli-Silvestre}. Using this extension, the non-singular case $\beta'\ge \delta>0$  was studied in~\cite{Athanasopoulos-Caffarelli}. If the equation is neither degenerate, $0< \delta\le\beta'\le C<\infty$,  one gets extra regularity; see~\cite{VPQR}. As for the singular case, the only precedent is~\cite{Kim-Lee}, where the authors consider the nonlinearity $\beta(s)=s^{1/m}$, $m\in(\frac{N-\sigma}{N+\sigma},1)$, $\sigma<N$, again by means of the extension, and prove H\"older regularity for nonnegative solutions for the problem posed in a bounded domain. The fact that the solution has a sign is crucial in their proof.

For general kernels, which do not admit an extension, a different point of view is needed. In the linear case two approaches have been used, one close to De Giorgi's ideas, see~\cite{Caffarelli-Chan-Vasseur},  and the other one more related to Moser's ones, see for instance~\cite{Felsinger-Kassmann}. We treat the degenerate case in~\cite{PQR} following the approach from~\cite{Caffarelli-Chan-Vasseur}. Let us remark that the  technique that we will use in the singular case works also in the degenerate case when the degeneracy is at most algebraic; see Theorem~\ref{thm.degenerate}.

\textsc{Outline of the regularity proof. }
The proof of regularity follows some ideas of the method introduced in the fifties of the last century by E. De Giorgi \cite{deGiorgi} to deal with elliptic equations. This approach, based on the control of the oscillation of the solution in a family of nested space-time cylinders, has been successfully applied with modifications to treat nonlinear local parabolic problems, see for instance~\cite{diBenedetto}, or linear nonlocal parabolic problems~\cite{Caffarelli-Chan-Vasseur,Caffarelli-Vasseur2}.

In order to take care of the nonlocal character of the operator we use ideas from~\cite{Caffarelli-Chan-Vasseur}. On the other hand, to deal with the nonlinearity we should look at~\cite{PQR}, where we considered the case of degenerate equations \eqref{eq:elefi}, where $\varphi'(0)=1/\beta'(0)=0$.
The main technical novelty with respect to that paper is that, instead of the quadratic energies that were used there, which coincide with the ones which are adequate to treat linear problems, here we need to use a \lq\lq nonlinear'' energy adapted to $\beta$, since in our case $\beta'(0)=0$. When $\mathcal{L}=(-\Delta)^{\sigma/2}$ and $\beta$ is a power, this energy coincides with the one used in~\cite{Kim-Lee}. Let us notice however that our treatment of the energy differs from the one therein, which is what allows us to consider sign changing solutions. In the local context such energies were introduced in~\cite{diBenedetto-Kwong}.

The first step in the regularity argument is to obtain a  De Giorgi type oscillation reduction lemma: if $u$ is mostly below a reference level in space-time measure in some parabolic cylinder, then the supremum goes down if we restrict to the half cylinder. Analogously, if $u$ is mostly above the reference level in space-time measure in some parabolic cylinder, then the infimum goes up if we restrict to the half cylinder.

In order to simplify the computations we set the reference level at 0 and the size at~1 and work with normalized cylinders. Therefore, we have to deal with solutions of scaled versions of equation~\eqref{eq:main}. In these scaled versions the singularity is translated to some unknown point. However,  we are able to obtain  an energy estimate that \emph{does not depend} on the location of the singularity.  This is the  energy inequality corresponding to the one obtained in the original work by De Giorgi which controls the $L^2$--norm of the gradient in terms of the size of the solution. Thus it is a kind of reversed Sobolev inequality. This is proved in Section \ref{sect-Energy}.

A second De Giorgi type lemma will tell us what happens when the solution is neither mostly positive nor mostly negative  in space-time measure. We will prove that some mass is lost between successive intermediate energy levels in $(0,1)$ to be defined,  a quantitative version of the fact that a function with a jump discontinuity cannot be in the energy space. Since we are away from the singularity on one of the sides of the reference level, this is a result of linear nature, which we borrow from~\cite{Caffarelli-Chan-Vasseur}.
Both De Giorgi type lemmas are included in Section \ref{sect-DeGiorgi}.

This is enough to prove the oscillation reduction result in Section~\ref{sect-Oscillation-Reduction}. The proof works as follows. Assume, without loss of generality, that the singular point (where the diffusion coefficient $1/\beta'$ is infinity) lies below 0, otherwise we work with $-u$. Then we will prove that eventually it will be possible to apply the first De Giorgi type lemma to one of the intermediate energy levels. Indeed, if this were not the case, we could apply the second De Giorgi type lemma to show that some fixed amount of energy would be lost between two successive  energy levels, which would lead to a contradiction after a finite number of steps.

From this we get next the continuity of the solution,   Theorem~\ref{th:regularity}, by means of scaling arguments. To prove the H\"older  regularity stated in Theorem~\ref{th:holder.continuity} we have to consider separately points where $u$ vanishes and points where $u$ is different from zero, since  the constants in the energy inequality degenerate for the rescaled problems when approaching a point of singularity. This is done in Section~\ref{sect-Regularity}.


\section{Energy inequality}

\label{sect-Energy} \setcounter{equation}{0}

We obtain in this section an energy estimate for the solutions of the rescaled versions of equation \eqref{eq:main} mentioned in the previous section. We want to deal with levels of the solution close to the singularity point $s_0$ and also levels far away from it. To treat them in an unified way we introduce the condition
\begin{equation}\label{eq:cond.beta3}\tag{${\rm H}_3$}
\beta\in C^1(\mathbb{R}),\qquad
\mathfrak{m}|s-s_0|^{\ell-1}\le \beta'(s)\le \mathcal{M} \quad \mbox{for } |s|\le 2,
\end{equation}
for some $s_0\in\mathbb{R}$, $\mathfrak{m},\,\mathcal{M}>0$ and $\ell\ge1$. In particular we will use $\ell>1$ close to the singularity, and $\ell=1$ away from it.

To get the energy inequality we will consider $\zeta=(u-\psi)_+$ as test function in \eqref{eq:main.weak}, for some space-dependent, nonnegative function $\psi$. Though $\zeta$ is not regular enough, in~\cite{PQR} it is shown, by means of some Steklov averages, that   the functional
$$
\mathcal{B}_{\psi}(u)=\int_{0}^{(u-\psi)_{+}}\beta^{\prime}(s+\psi)s\,ds
$$
satisfies
\begin{equation}
\label{weak-u}\left. \int_{\mathbb{R}^{N}}\mathcal{B}_{\psi}(u(x,t))\,dx\right| _{t_{1}%
}^{t_{2}}+\int_{t_{1}}^{t_{2}}\mathcal{E}_J(u,(u-\psi)_{+})(t)\,dt=0.
\end{equation}
If $\psi$ is a Lipschitz function satisfying  $\int_{\{|x-y|>1\}}|\psi(x)-\psi(y)|J(x,y)\,dy<C<\infty$ for every $x\in\mathbb{R}^N$, the second term is estimated in \cite{Caffarelli-Chan-Vasseur} by
\begin{equation}\label{eq:E(u,u)}
  \mathcal{E}_J(u,(u-\psi)_+)\ge \mathcal{E}_J((u-\psi)_+,(u-\psi)_+)-C\big(\|(u-\psi)_+\|_1+|\{u>\psi\}|\big).
\end{equation}
Observe  that if $\psi$ is a constant, the estimate would be easier
$$
\mathcal{E}_J(u,(u-\psi)_+)\ge \mathcal{E}_J((u-\psi)_+,(u-\psi)_+),
$$
use for instance Stroock-Varopoulos inequality, see~\cite{Varopoulos,Brandle-dePablo}.
The introduction of an unbounded barrier function like $\psi$ is needed due to the nonlocal character of the equation. Indeed, $\zeta\in C^\infty_{\rm c}(\mathbb{R}^N)$ does not imply that $\mathcal{L}\zeta$ is compactly supported. Hence, though we are considering bounded solutions, the successive versions obtained in the scaling procedure do not have a uniform bound.  The idea to \emph{localize} these nonlocal problems is to impose the growth condition at infinity $u(x,t)\le \psi^{1/2}(x)$, which  is kept under the scaling. The set of points where $u>\psi$ is uniformly localized for solutions satisfying that condition, which implies a uniform bound for all of them, see \cite{Caffarelli-Chan-Vasseur}.

We now look at the first term in \eqref{weak-u}. The following calculus result is crucial in what follows.
\begin{Proposition}\label{prop-area} If $\beta$ satisfies \eqref{eq:cond.beta3} there exists a positive constant $c_\ell$ depending only on $\ell$ such that for every $0\le k<w\le2$ it holds
\begin{equation*}
\label{estimate-B}c_\ell \mathfrak{m}(w-k)^{\ell+1}\le
\int_{0}^{w-k}\beta'(s+k)s\,ds\le \frac{\mathcal{M}}2 (w-k)^2.
\end{equation*}
\end{Proposition}
\noindent\emph{Proof. }
The upper estimate is trivial.
In order to estimate the integral from below we define the quantity
$$
r=\max\{0,\,\min\{w-k,\,s_0-k\}\,\}.
$$
Using \eqref{eq:cond.beta3} we have to estimate the function $F(s)=|s+k-s_0|^{\ell-1}s$. It holds
$$
F(s)\ge\begin{cases}
2^{1-\ell}r^{\ell-1}s,&\mbox{ for } 0< s\le r/2, \\
(s-r)^{\ell},&\mbox{ for } r\le s< w-k.
\end{cases}
$$
Thus
$$
\begin{array}{rl}
\displaystyle\int_0^{w-k} F(s)\,ds&\displaystyle\ge\int_0^{r/2}2^{1-\ell}r^{\ell-1}s\,ds+\int_r^{w-k} (s-r)^{\ell}\,ds \\ [3mm]
&\displaystyle=\frac1{2^{\ell+2}}r^{\ell+1}+\frac1{\ell+1}(w-k-r)^{\ell+1}\ge \frac1{2^{2\ell+2}} (w-k)^{\ell+1}.
\end{array}$$
\qed

As a consequence, we get for the energy functional the estimate
\begin{equation}
\label{estimate-B2}c_\ell \mathfrak{m}(u-k)_+^{\ell+1}\le
B_\psi(u)\le \frac{\mathcal{M}}2 (u-k)_+^2.
\end{equation}

Observe that $\|f\|_2+\mathcal{E}_J^{1/2}(f,f)$ is a norm equivalent to the standard Sobolev norm $\|f\|_2+\|(-\Delta)^{\sigma/4}f\|_2$ in $H^{\sigma/2}(\mathbb{R}^N)$. We recall that in our definition of weak solution we do not require  $u\in L^2(\mathbb{R}^{N})$, but $\beta(u)\in L^1(\mathbb{R}^{N})$. Nevertheless, the localization performed by using the barrier functions gives $(u-\psi)_+\in  L^q(\mathbb{R}^{N})$ for every $q\ge1$.  Therefore, putting together \eqref{weak-u}, \eqref{eq:E(u,u)} and \eqref{estimate-B2}, we obtain the desired energy inequality.
\begin{Lemma}\label{lem:newenergy} Assuming $J$ and $\beta$ satisfy, respectively, \eqref{eq:kernel} and \eqref{eq:cond.beta3},  there exist positive constants $C_1,\,C_2$ such that if $u$ is a bounded weak solution to \eqref{eq:main} satisfying $|\{u(\cdot,t)>\psi(x)\}|<\infty$, where $\psi$ is a barrier as before, then
\begin{equation}\label{eq:newenergy}
\begin{array}
{rl}%
\displaystyle \mathfrak{m}\|(u-\psi)_+(\cdot,t_2)\|_{\ell+1}^{\ell+1}&\displaystyle+ \int_{t_{1}%
}^{t_2}\|(-\Delta)^{\sigma/4}((u-\psi)_{+})(\cdot,t)\|_2^2\,dt\\[10pt]%
&\displaystyle\le  C_1\mathcal{M}\|(u-\psi)_+(\cdot,t_1)\|_2^2+C_2\int_{t_{1}}^{t_2} z(t)\,dt.
\end{array}
\end{equation}
where
$$
z(t)=\|(u-\psi)_+(\cdot,t)\|_2^2+\|(u-\psi)_+(\cdot,t)\|_1+|\{u(\cdot,t)>\psi(x)\}|.
$$
\end{Lemma}
This  inequality allows to control a nonlinear energy of the truncated function $u_\psi=(u-\psi)_+$ in terms of its $L^2$ norm.


\section{De Giorgi type lemmas}

\label{sect-DeGiorgi} \setcounter{equation}{0}

The next step is to obtain a first De Giorgi type oscillation reduction lemma: if $u$ is mostly negative in space-time measure in some parabolic cylinder, then the supremum goes down if we restrict to the half cylinder.

\medskip

\noindent\emph{Notation. }
$\Gamma_{R,a}=\{|x|<R,\,-a\le t\le 0\}$.

\begin{Lemma}
\label{lem:first.DG.lemma}
Assume \eqref{eq:cond.beta3} holds. There is a constant $\delta\in(0,1)$ such that,  if $u:\mathbb{R}^N\times[-2,0]\to \mathbb{R}$ is a  weak solution to equation~\eqref{eq:main} satisfying, for some $0<a\le1$,
\begin{align}
\label{eq:primerapsi}
&u(x,t)\le 1+(|x|^{\sigma/4}-1)_{+}\quad\text{in } \mathbb{R}^N\times[-2,0],\\[10pt]
\label{eq:critical.condition}
&|\{u>0\}\cap \Gamma_{2,2a}|\le \delta a^{2(1+N/\sigma)},
\end{align}
then
$$
u(x,t)\le 1/2 \quad\text{if }(x,t)\in \Gamma_{1,a}.
$$
\end{Lemma}

\noindent\emph{Proof. }
Let $L_k=\frac12-\frac{1}{2^{k+1}}$, $t_k=-(1+\frac{1}{2^k})a$. We consider the sequence of barriers $\psi(x)=\psi_{L_k}(x)=L_k+(|x|^{\sigma/2}-1 )_+$ in~\eqref{eq:newenergy}, which satisfy the requirements needed to get \eqref{eq:E(u,u)}. Let $u_k(t)=(u-\psi_{L_k})_{+}(\cdot,t)$. We define the quantity
\begin{equation*}\label{eq:Uk}
U_{k}=\sup_{t_{k}\le t\le 0}\|u_{k}(t)\|_{\ell+1}^{\ell+1}+ \int_{t_{k}}^{0}\|(-\Delta)^{\sigma/4}u_k(t)\|_2^2\,dt,
\end{equation*}
corresponding to the different nonlinear energy levels of the truncated functions $u_k$.
Observe that $(|x|^{\sigma/2}-1 )_+<u\le 1+(|x|^{\sigma/4}-1)_{+}$ implies $|x|<(\frac{1+\sqrt5}2)^{4/\sigma}$ and $0<u<\frac{1+\sqrt5}2$, so that we can use condition \eqref{eq:cond.beta3} and Lemma \ref{lem:newenergy}. First, the energy estimate~\eqref{eq:newenergy} implies, for $k\ge1$ and $t_{k-1}<s<t_k$,
$$
U_k\le  (1+\frac1{\mathfrak{m}})\left(C_1\mathcal{M}\|u_k(s)\|_2^2+C_2\int_{t_{k-1}}^{0}z_k(\tau)\,d\tau\right).
$$
Taking the mean in the interval $s\in[t_{k-1},t_k]$, we get
\begin{equation}\label{eq:Uk2}
\begin{array}{rl}
U_k&\displaystyle\le \frac{1+1/\mathfrak{m}}{t_k-t_{k-1}}\int_{t_{k-1}}^{t_k}\left(C_1\mathcal{M}\|u_k(s)\|_2^2+
C_2\int_{t_{k-1}}^{0}z_k(\tau)\,d\tau\right)\,ds \\[10pt]%
&\displaystyle\le \frac{C\kappa2^{k}}a\int_{t_{k-1}}^{0}\Big(\|u_k(s)\|_2^2+ \|u_k(s)\|_1+|\{u_k(s)>0\}|\Big)\,ds,
\end{array}
\end{equation}
where $\kappa=(\mathcal{M}+1)(1+1/\mathfrak{m})$.

Now, since $L_k=L_{k-1}+\frac 1{2^{k+1}}$, then $u_k>0$ implies $u_{k-1}>\frac 1{2^{k+1}}$, which in turn gives the Chebyshev type inequality
$$
\int_{\mathbb{R}^N}u_k^r\le 2^{(k+1)(q-r)}\int_{\mathbb{R}^N}u_{k-1}^q
$$
for every $q>r$. Thus, for some $q\ge2$ to be chosen we get, using this inequality with $r=0,\,1,\,2$, that \eqref{eq:Uk2} reduces to
\begin{equation}\label{eq:Uk3}
U_k\le \frac{C\kappa2^{(q+1)k}}a\int_{t_{k-1}}^{0}\|u_{k-1}(t)\|_q^qdt.
\end{equation}

To link this estimate with $U_{k-1}$ we use Hardy-Littlewood-Sobolev inequality, so $N>\sigma$ is required.  Using first interpolation we get, with $q=\frac{(\ell+1)\sigma}{N}+2$,
$$
\begin{array}
[c]{rl}%
\displaystyle\int_{t_{k-1}}^{0}\|u_{k-1}(t)\|_{q}^{q}\,dt & \displaystyle\le
\int_{t_{k-1}}^{0}\|u_{k-1}(t)\|_{\ell+1}^{(\ell+1)\sigma/N}\|u_{k-1}(t)\|_{\frac
{2N}{N-\sigma}}^{2}\,dt\\[4mm]
& \displaystyle\le\left( \sup_{t_{k-1}<t<0}\|u_{k-1}(t)\|_{\ell+1}^{\ell+1}\right)
^{\sigma/N} \int_{t_{k-1}}^{0}\|(-\Delta)^{\sigma/4}u_{k-1}(t)\|_{2}^2\,dt\\[4mm]
& \displaystyle\le C\left( \sup_{t_{k-1}<t<0}\|u_{k-1}(t)\|_{\ell+1}^{\ell+1}+ \int_{t_{k-1}%
}^{0}\|(-\Delta)^{\sigma/4}u_{k-1}(t)\|_2^2\,dt\right)
^{1+\sigma/N}\\[6mm]
& \displaystyle\le C U_{k-1}^{1+\sigma/N}.
\end{array}
$$
In the case $N=1\le\sigma$ we use a Nash-Gagliardo-Nirenberg inequality, see~\cite{PQRV2}, to get the same estimate. We have thus arrived to the following nonlinear recurrence
\begin{equation*}\label{eq:ukuk-1}
U_k\le \frac{C\kappa2^{(q+1)k}}a U_{k-1}^{1+\sigma/N}.
\end{equation*}
Therefore there exists $\varepsilon>0$ such that if $ U_1<\varepsilon (a/\kappa)^{N/\sigma} $ then  $U_k\to0$ as $k\to\infty$; see \cite{diBenedetto}. This will give $u<\psi_{L_\infty}$, that is $u(x,t)<1/2$ for $|x|<1$, $-a\le t\le 0$.

The condition on $U_1$ is fulfilled, thanks to estimate~\eqref{eq:Uk3} in the particular case $k=1$ and $q=2$, if
\begin{equation}\label{eq:u0small}
\int_{-2}^0\int_{\mathbb{R}^N}(u-\psi_{L_0})_+^2<\varepsilon \left(\frac a\kappa\right)^{1+N/\sigma}.
\end{equation}
But this is not guaranteed by condition \eqref{eq:critical.condition} since $\{u>\psi_{L_0}\}\not\subset \Gamma_{2,2a}$,
and we need to apply a scaling argument.

Let $(x_0,t_0)\in \Gamma_{1,a}$ be arbitrary, and define for some large $\Lambda$ the function
\begin{equation*}
  \label{eq:uR}
  u_\Lambda(x,t)=u(x_0+\Lambda^{-1}x,t_0+\Lambda^{-\sigma}t).
\end{equation*}
This function solves  equation \eqref{eq:main} with
the new bilinear form $\mathcal{E}_{J_\Lambda}$  associated to the rescaled kernel
$$
J_\Lambda(x,y)=\Lambda^{-(N+\sigma)}J(x_0+\Lambda^{-1}x,x_0+\Lambda^{-1}y),
$$
which satisfies again hypothesis~\eqref{eq:kernel} with the same constant whenever $\Lambda\ge1$.
On the other hand,  the following property is proved in \cite{Caffarelli-Chan-Vasseur},
$$
1+(|x_0+\Lambda^{-1}x|^{\sigma/2}-1)_{+}\le (|x|^{\sigma/4}-1)_{+},\quad \text{for every } |x_0|<1,\;|x|>\Lambda,
$$
if $\Lambda$ is chosen large enough depending only on $\sigma$.
Thus, since $u$ satisfies condition \eqref{eq:primerapsi} we have that $u_\Lambda$ satisfies $u_\Lambda(x)\le\psi_{L_0}(x)$ for every $|x|>\Lambda$, and therefore
$$
\begin{array}{l}
\displaystyle\int_{-2}^0\int_{\mathbb{R}^N}\big(u_\Lambda(x,t)-\psi_{L_0}(x))_+\big)_+^2\,dxdt
\le\int_{-2}^0\int_{\{u_\Lambda(\cdot,t)>0\}\cap B_\Lambda}u^2_\Lambda(x,t)\,dxdt \\ [10pt]
\displaystyle\qquad\le \Lambda^{N+\sigma}\int_{t_0-2\Lambda^{-\sigma}}^{t_0}\int_{\{u(\cdot,\tau)>0\}\cap B_{1}(x_0)}u^2(z,\tau)\,dzd\tau \\ [12pt]
\displaystyle\qquad
\le \Lambda^{N+\sigma}\int_{-2a}^0\int_{\{u(\cdot,\tau)>0\}\cap B_2}u^2(z,\tau)\,dzd\tau \le4\Lambda^{N+\sigma}|\{u>0\}\cap B_2\times[-2a,0]|.
\end{array}
$$
We have used that  $|z-x_0|<1$ implies $|z|<2$ and $t_0-2\Lambda^{-\sigma}>-2a$ if $\Lambda\ge(2/a)^{1/\sigma}$. Choosing $\Lambda=(2/a)^{1/\sigma}$ and  $\delta=c\varepsilon\kappa^{-(1+N/\sigma)}$ in \eqref{eq:critical.condition}, we get  that $u_\Lambda$ satisfies \eqref{eq:u0small} and thus $u_\Lambda<1/2$ in $\Gamma_{1,a}$. In particular $u_\Lambda(0,0)<1/2$, which means $u(x_0,t_0)<1/2$.

\qed

\begin{Remark}\label{rem-rem}
The appearance of the constant $\kappa$ is an effect of the nonlinearity. Observe that this constant is large, and therefore $\delta$ is small, when $\mathfrak{m}$ is small (even if $\mathcal{M}$ were also small). This will be the case for the rescaled problems considered in the proof of H\"older regularity close to a singular point, where the rescaled values of $\mathfrak{m}$ and $\mathcal{M}$ go to zero.
\end{Remark}

Applying this lemma to $-u$ we get also that
if $u$ is mostly positive in space-time measure in the cylinder $\Gamma_{2,2a}$, then the infimum goes up in $\Gamma_{1,a}$. Observe that $-u$ solves problem~\eqref{eq:main} with $\beta$ replaced by $\widetilde\beta(s)=-\beta(-s)$. In both cases we have reduced the oscillation of the solution in the half cylinder.

To proceed with the regularity proof we need to analyze what happens when the solution is neither mostly positive nor negative, in the sense of Lemma~\ref{lem:first.DG.lemma}, in space-time measure.

To this aim we will use De Giorgi's idea of loss of mass at intermediate levels.
This comes from a result in  \cite{Caffarelli-Chan-Vasseur} on the intermediate values.
The key idea is to impose conditions on the nonlinearity guaranteeing that the equation is not singular at the intermediate values. Hence we are in the linear setting studied in~\cite{Caffarelli-Chan-Vasseur}, and the proof goes as there. The result is written in terms of  a cut-off function $F$, continuous radially nonincreasing such that $F\equiv1$ for $0\le |x|\le1$, $F\equiv0$ for $|x|\ge2$, and the functions
$$
\psi_\lambda(x)=((|x|-\lambda^{-4/\sigma})_+^{\sigma/4}-1)_+,
$$
used to control the growth at infinity. Observe that  $\psi_\lambda\equiv0$ for $|x|<c(\lambda)$, with $c(\lambda)>10$ if  $\lambda<1/3$.

\begin{Lemma}
\label{lem:second.DG.lemma} Assume $0<C_1\le \beta'(s)\le C_2$ for every $1/2\le s\le2$. For every $\nu,\mu>0$ small enough, $0<a<1$, there exist $\gamma>0$ and $\bar\lambda\in(0,1/3)$ such that for any $\lambda\in (0,\bar\lambda)$, and any weak solution $u:\mathbb{R}^N\times[-2,0]\to\mathbb{R}$ to~\eqref{eq:main} satisfying
$$
u(x,t)\le 1+\psi_\lambda(x)\quad\text{on }\mathbb{R}^N\times[-2,0],\qquad |\{u<0\}\cap(B_1\times (-2,-2a))|\ge\mu,
$$
we have the following implication: If
$$
|\{u>1-\lambda^2 F\}\cap \Gamma_{2,2a}|\ge\nu,
$$
then
$$
|\{ 1-F<u<1-\lambda^2 F\}\cap \Gamma_{2,2}|\ge\gamma.
$$
\end{Lemma}


\section{Oscillation reduction}

\label{sect-Oscillation-Reduction} \setcounter{equation}{0}

 In order  to apply Lemma~\ref{lem:second.DG.lemma} to the different energy levels,  we need a more restrictive control of the behaviour of the solution at infinity, which is given in terms of
$$
H_{\lambda}(x)=\big((|x|-\lambda^{-4/\sigma})_+^{\sigma\lambda^2/4}-1\big)_{+},\qquad
\lambda>0.
$$

\begin{Lemma}\label{lema-oscil}
Assume  that \eqref{eq:cond.beta3} holds, and let $\bar\lambda$ be as in Lemma~\ref{lem:second.DG.lemma}. There exist constants $\tau,\theta\in(0,1)$ such that if
$u:\mathbb{R}^N\times[-2,0]\to \mathbb{R}$ is a  weak solution to equation~\eqref{eq:main}   that satisfies, for $\lambda\in(0,\bar\lambda)$
\begin{equation}
\label{eq:condition.oscillation.reduction}
|u(x,t)|\le1+ H_\lambda(x)\quad\text{on }\mathbb{R}^N\times[-2,0],
\end{equation}
then
$$
\sup_{\Gamma_{1,\tau}} u-\inf_{\Gamma_{1,\tau}} u\le2-\theta.
$$
\end{Lemma}

\noindent\emph{Proof. } If $u$ or $-u$ satisfy estimate \eqref{eq:critical.condition} for some $0<a\le1$, we are done since then Lemma \ref{lem:first.DG.lemma} gives  the result with $\theta=1/2$ and $\tau=a$ (observe that condition~\eqref{eq:condition.oscillation.reduction} implies condition~\eqref{eq:primerapsi}).  Now assume for instance $s_0\le0$, so we have that the nonlinearity satisfies the hypotheses of Lemma~\ref{lem:second.DG.lemma}. Then, since  $-u$ does not satisfy \eqref{eq:critical.condition} for any $0<a\le1$, taking $a=1$  we  have
$|\{u<0\}\cap \Gamma_{2,2}|>\delta$.
Thus, $|\{u<0\}\cap B_2\times (-2,-b)|\ge\delta/2>0$ for some $b\in (0,\frac{\delta}{2|B_2|})$, and thus we can apply Lemma~\ref{lem:second.DG.lemma}.

We consider now the sequence of rescaled functions
$$
u_{k+1}=\frac{u_k-(1-\lambda^2)}{\lambda^2},\qquad u_0=u.
$$
Then $u_k$ satisfies
$$
\partial_t \beta_k(u_k)+\mathcal{L} u_k=0,
$$
with a nonlinearity $\beta_{k}$ given iteratively by
$$
\beta_{k+1}(s)=\frac{1}{\lambda^2}\beta_k(\lambda^2 s+ 1-\lambda^2),\qquad \beta_0=\beta,
$$
always with the same operator $\mathcal{L}$. We will prove that for each $k\ge1$ we can apply either Lemma~\ref{lem:first.DG.lemma} or Lemma~\ref{lem:second.DG.lemma}. Repeated application of Lemma~\ref{lem:second.DG.lemma} will give that in fact Lemma~\ref{lem:first.DG.lemma} can be applied after a finite number of steps. Hence we will be done.

The key point is that $\beta_{k+1}'(u_{k+1})=\beta_k'(u_{k})$, and $u_{k+1}>0$ implies $u_k>1-\lambda^2>0$. Thus \eqref{eq:cond.beta3} holds with some singularity point $s_k=1+(s_0-1)\lambda^{-2k}<0$. Also, if $u_{k+1}\in[1/2,2]$ then $u_{k}\in[1/2,2]$. We have in this way that the hypotheses on the nonlinearity of Lemmas \ref{lem:first.DG.lemma} and \ref{lem:second.DG.lemma} hold. In fact we may put $\ell=1$ in Lemma \ref{lem:first.DG.lemma}, since $s_0\le0$ and $u\ge1/2$. On the other hand, $H_\lambda(x)\le\lambda^{2} \psi_\lambda(x)$, and since
$$
u_{k+1}(x,t)\le 1+\frac{u_k(x,t)}{\lambda^2},
$$
we get by induction that $u_k(x,t)\le1+\psi_\lambda(x)$.

We now have
$$
|\{u_{k+1}<0\}\cap B_2\times(-2,-b)|\ge |\{u_k<0\}\cap B_2\times(-2,-b)| \ge\frac\delta2.
$$
Then, applying Lemma~\ref{lem:second.DG.lemma} with some $\nu>0$ to be chosen, we get that there exists $\gamma>0$ such that
$$
\begin{array}{l}
|\{u_{k+1}>1-\lambda^2 F\}\cap \Gamma_{2,b}|\\
\qquad\qquad
=|\{u_{k+1}>1-F\}\cap \Gamma_{2,b}|-|\{1-\lambda^2F>u_{k+1}>1-F\}\cap \Gamma_{2,b}|\\
\qquad\qquad\le|\{u_k>1-\lambda^2 F\}\cap \Gamma_{2,b}|-\gamma\le |\{u>1-\lambda^2 F\}\cap \Gamma_{2,b}|-k\gamma,
\end{array}
$$
and we arrive to a contradiction if  $\overline k\ge|\Gamma_{2,2}|/\gamma$. Therefore, and here we consider the particular value of $\nu= cb^{2+N/\sigma}$, there exists some $k_*< \overline k$ for which
$$
|\{u_{k^*}>0\}\cap \Gamma_{2,b}|\le c\delta a^{1+N/\sigma},
$$
as needed in Lemma~\ref{lem:first.DG.lemma} with $a=b/2$. This means
$$
u_{k_*}\le \frac12\quad\text{in }\Gamma_{1,b/2}.
$$
Going back to the original variables we get in that set that
$$
-1\le u = 1+\lambda^{2k_*}(u_{k_*}-1)\le 1-\theta, \qquad\theta=\lambda^{2|\Gamma_{2,2}|/\gamma}/2,
$$
so the oscillation in $\Gamma_{1,\tau}$ is less than $2-\theta$ with $\tau=b/2$. A final comment on the reduction constant $\theta$: it is small when $\gamma$ is small, that is, when $\delta$ is small.
\qed


\section{Regularity}

\label{sect-Regularity} \setcounter{equation}{0}

Lemma~\ref{lema-oscil}  shows that the oscillation of $u$ in $\Gamma_{c(\lambda),2}$ is reduced in $\Gamma_{1,\tau}$ by a factor $\theta^*=1-\theta/2$. From this we get next the regularity stated in  Theorems~\ref{th:regularity} and ~\ref{th:holder.continuity}.

\noindent\emph{Proof of Theorem~\ref{th:regularity}} We want to prove regularity of the solution at an arbitrary given point. The first step is a traslation that moves that point to the origin. We may also assume that the $L^\infty$--norm of the solution is one. More precisely, let $(x_0,t_0)\in Q$ and $\tau_0=\inf\{1,t_0/3\}$, and put $A=\|u(\cdot,0)\|_\infty$. Then
$$
u_0(x,t)=\frac1Au(x_0+\tau_0^{1/\sigma}x,t_0+\tau_0t)
$$
is a solution to the equation
$$
\displaystyle \int_a^b\int_{\mathbb{R}^{N}%
}\beta_0(u_0)\partial_{t} \zeta-\int_a^b\mathcal{E}_{J_0}(u_0,\zeta
)=0$$
for every $-2<a<b<0$, where $\beta_0(s)=\frac{1}{A}\beta(As)$, and
$$
J_0(x,y)=\tau_0^{\frac{N+\sigma}\sigma}J(x_0+\tau_0^{1/\sigma}x,x_0+\tau_0^{1/\sigma}y).
$$
The function $\beta_0$ and the kernel $J_0$ satisfy the
same hypotheses as $\beta$ and $J$.

Let now $Q_k=\Gamma_{R^{-k}, R^{-k\sigma}}$ for every $k\ge0$ and for some $R>1$ large enough to be determined later. We will show that  the semi-oscillation of $u$ in $Q_{k-1}$,
$$
\displaystyle\varpi_k=\frac{\sup_{Q_{k-1}} u-\inf_{Q_{k-1}} u}{2},
$$
goes to 0 as $k\to\infty$, which yields the result. To this aim, we assume from now on that $\varpi_{k}\ge \varsigma>0$, and we will arrive to a contradiction.

Given $k\ge1$, we define
\begin{equation}\label{iter1}
\displaystyle u_{k}(x,t)=\frac{u(R^{-k}x,R^{-\sigma k}t)-\mu_{k}}{\varpi_{k}},\quad \mu_k=\frac{\sup_{Q_{k-1}} u+\inf_{Q_{k-1}} u}{2}.
\end{equation}
The functions $u_k$ satisfy the equation
$$
\displaystyle \int_a^b\int_{\mathbb{R}^{N}%
}\beta_k(u_k)\partial_{t} \zeta-\int_a^b\mathcal{E}_{J_k}(u_k,\zeta
)=0$$
where
$$
\beta_{k}(s)=\dfrac{\beta\left(\varpi_{k} s+\mu_{k}\right)}{\varpi_{k}},\quad
J_{k}(x,y)=R^{-(N+\sigma)k}J_0(R^{-k}x,R^{-k}y).
$$
Observe that $J_k$  satisfies again~\eqref{eq:kernel}. Notice also that, using~\eqref{eq:cond.beta1}, we have
$$
\beta_k'(s)\ge \mathfrak{m}_0\varpi_k|s-s_k|^{p-1}, \quad\text{where }s_k=-\mu_k/\varpi_{k}.
$$
We have in particular that all the nonlinearities $\beta_k$ satisfy~\eqref{eq:cond.beta3} with the same constants $\mathfrak{m}=\mathfrak{m}_0\varsigma$, $\mathcal{M}=\mathcal{M}_0$ and $\ell=p$.

Now we check that all the functions $u_k$ satisfy condition~\eqref{eq:condition.oscillation.reduction}. This is clear for $|x|\le R$, since  $|u_k(x,t)|\le 1$. For $|x|>R$ we have
$$
|u_k(x,t)|\le \frac{2}{\varpi_k}\le \frac2\zeta,
$$
so it is enough to take $R$ large such that $H_\lambda(R)\ge\frac{2-\zeta}\zeta$.

Applying Lemma \ref{lema-oscil} we get
$$
\frac{\sup_{\Gamma_{1,\tau}} u_k-\inf_{\Gamma_{1,\tau}} u_k}{2}\le \theta^*.
$$
Therefore, taking $R^\sigma>1/\tau$ we conclude
$$
\begin{array}{rl}
\displaystyle\varpi_{k+2}&\displaystyle=\frac{\sup_{\Gamma_{R^{-(k+1)},R^{-(k+1)\sigma}}} u-\inf_{\Gamma_{R^{-(k+1)},R^{-(k+1)\sigma}}} u}{2}
\\ [4mm] &\displaystyle\le
\frac{\sup_{\Gamma_{R^{-k},R^{-k\sigma}\tau}} u-\inf_{\Gamma_{R^{-k},R^{-k\sigma}\tau}} u}{2}
\\ [4mm] &\displaystyle\le \frac{\sup_{\Gamma_{1,\tau}} u_k-\inf_{\Gamma_{1,\tau}} u_k}{2}\,\varpi_k\le
\theta^*\varpi_k.
\end{array}
$$
This is the desired contradiction.
\qed

In order to prove H\"older regularity we assume that the nonlinearity $\beta$ does not oscillate too much at the origin, that is, it behaves like a power times a slowly varying function.

\noindent\emph{Proof of Theorem~\ref{th:holder.continuity}. } As before we assume $\|u\|_\infty=1$ and show regularity at the origin.

\noindent\textsc{H\"older regularity at nonsigular points. } Suppose $u(0,0)>0$, the case $u(0,0)<0$ being similar.
We define  a
sequence of functions similar to~\eqref{iter1}, though in this case we divide by an estimate of the oscillation, instead of the oscillation itself. We take profit of the continuity that we have just proved.

Let $k_0\ge1$ be such that $\varpi_{k_0}s+\mu_{k_0}\ge\frac{u(0,0)}2>0$ for every $|s|\le2$, and take
\begin{equation*}\label{iter3}
\displaystyle u_{k}(x,t)=\frac{u(R^{-k}x,R^{-\sigma k}t)-\mu_{k}}{\nu_{k}},\qquad \nu_{k}=\begin{cases}
  \varpi_{k}, &1\le k\le k_0,\\ \varpi_{k_0}\theta_*^{k-k_0},&k>k_0,
\end{cases}
\end{equation*}
$\varpi_k$ and $\mu_k$ as before, $R$ large to be fixed later.
The functions $u_k$  satisfy  the equation
$$
\displaystyle \int_a^b\int_{\mathbb{R}^{N}%
}\beta_k(u_k)\partial_{t} \zeta-\int_a^b\mathcal{E}_{J_k}(u_k,\zeta
)=0$$
where the new nonlinearity is
$$
\beta_{k}(s)=\dfrac{\beta\left(\nu_ks+\mu_k\right)}{\nu_k},
$$
and the kernel $J_k$ is as before.  We must check that both $u_k$ and $\beta_k$ satisfy the hypotheses of
Lemma~\ref{lema-oscil}, and we start with $k=k_0$.

As to $\beta_k$ it is clear that it satisfies \eqref{eq:cond.beta3}  with $\ell=1$ and constants $\mathfrak{m}=\mathfrak{m}_0(\frac{u(0,0)}2)^{p-1}$, $\,\mathcal{M}=\mathcal{M}_0$, since $\nu_ks+\mu_k\ge\frac{u(0,0)}2$  for every $|s|\le2$, $k\ge k_0$.
Now we look at the sequence $u_k$. When $|x|<R$, $-2<t<0$ we have $|u_k|\le 1\le 1+H_\lambda(x)$ if $R^\sigma>2/\tau$. Let then be $|x|>R$. We have
$$
u_{k_0}(x,t)\le\frac{1+H_\lambda(x/R^{k_0})-\mu_{k_0}}{\varpi_{k_0}}\le 1+H_\lambda(x),
$$
provided $R$ is large, and an analogous estimate from below. Now assume by induction that
$$
u_{k}(x,t)\le 1+H_\lambda(x),
$$
for some $k\ge k_0$. We get, for $|x|>R$,
$$
\begin{array}{rl}
\displaystyle u_{k+1}(x,t)&\displaystyle\le\frac{u(R^{-(k+1)}x,R^{-\sigma(k+1)}t)-\mu_{k+1}}{\nu_{k+1}} \\ [3mm]
&\displaystyle\le\frac{\mu_k-\mu_{k+1}+\nu_k(1+H_\lambda(x/R))}{\nu_{k+1}} \\ [3mm]
&\displaystyle\le \frac{2(1+H_\lambda(x/R))}{\theta^*}\le 1+H_\lambda(x),
\end{array}
$$
again if $R$ is large enough.

We conclude an oscillation estimate of order
$(\theta^*)^k$ for $u$ in $Q_k$. This gives  H\"older regularity at points where the equation is nonsingular. Notice  that $\theta^*$ depends on $\mathfrak{m}$, which degenerates as $u$ approaches the value zero, see Remark~\ref{rem-rem} and the comment at the end of the proof of Lemma~\ref{lema-oscil}.

\noindent\textsc{H\"older regularity at singular points. }
Let now $u(0,0)=0$.
We assume here that \eqref{eq:cond.beta2} holds true and consider the
sequence of functions defined by means of a recurrence that takes into account the nonlinearity, and the singularity of $1/\beta'$ at zero. Also, since in that case the absolute vale  is controlled by the oscillation we may avoid substracting the mean. We define, for some $0<\gamma<1$ to be chosen,
$$
\displaystyle u_{k+1}(x,t)=\frac1{\gamma}\,u_k\left(\dfrac xR,\dfrac t{\gamma R^{\sigma}}\right),\qquad u_0=u.
$$
The corresponding rescaled nonlinearity in the problem satisfied by $u_k$ turns to be
$$
\beta_{k}(s)=\dfrac{\beta(\gamma^ks)}{\beta(\gamma^k)}.
$$
Our  goal is to prove that the oscillation of $u$ in each cube as before is $\varpi_k\le c\gamma^k$ for $k\ge 1$, thus implying H\"older regularity.   To estimate the oscillation we check again that the pairs $(u_k,\beta_k)$ fulfill the conditions of Lemma~\ref{lema-oscil}. Observe first that by \eqref{eq:cond.beta2}
$$
\beta_k'(s)=\frac{\gamma^k\beta'(\gamma^ks)}{\beta(\gamma^k)}\sim
\frac{|s|^{q-1}\gamma^{kq}h(\gamma^ks)}{\int_0^{\gamma^k} r^{q-1}h(r)\,dr}\sim |s|^{q-1},
$$
where we have used the properties of slowly varying functions, \cite{Karamata}. We thus get  \eqref{eq:cond.beta3} with $\ell=q$ and $s_0=0$ and $\mathfrak{m}$ bounded away from zero. On the other hand, by induction applying Lemma~\ref{lema-oscil} to $u_k$ we know that $|u_k(x,t)|\le\theta^*$ for $|x|<1$, $-\tau<t<0$, where $\theta^*$ depends only on $\mathfrak{m}$ and not on $\gamma$. We therefore may put $\gamma=\theta^*$. Thus $|u_{k+1}(x,t)|\le1$ for $|x|<R$, $-2<t<0$ if we take $\tau \theta^*R^\sigma>2$. Outside the ball, $|x|>R$, we have
$$
|u_{k+1}(x,t)|\le\frac1{\theta^*}\left|\,u_k\left(\dfrac xR,\dfrac t{\theta^* R^{\sigma}}\right)\right|\le \frac{1+H_\lambda(x/R)}{\theta^*}\le1+H_\lambda(x)
$$
provided $R$ is large. The proof is complete.
\qed

\section{Degenerate equations}

\label{sect-Degenerate} \setcounter{equation}{0}

The approach that we have followed can be applied to degenerate equations, when the diffusion coefficient $1/\beta'(u)$ vanishes at zero, giving a simplified proof of the results in~\cite{PQR}. Let us remark however that the conditions on the nonlinearities there and here are not exactly the same ones. We thus arrive to the following result.
\begin{Theorem}
\label{thm.degenerate}
Let $\beta\in C(\mathbb{R})\cap C^1(\mathbb{R}\setminus\{0\})$ satisfy $\beta(0)=0$ and
$$
\beta'(s)\ge\mathfrak{m}_0\;\mbox{ for } s\in\mathbb{R},\quad\beta'(s)\le \mathcal{M}_0|s|^{p-1}\;\mbox{ for } |s|\le \|u\|_\infty,\quad \mathfrak{m}_0>0,\,p\in(0,1),
$$
then $u$ is continuous. If moreover $\beta$ satisfies \eqref{eq:cond.beta2}, then $u$ is H\"older continuous at every point.
\end{Theorem}

The main idea is that if a function $\varphi$ satisfies
  $$
  \mathfrak{m}\le \varphi'(s)\le \mathcal{M}|s-s_0|^{\ell-1}, \qquad \mbox{for every } |s|\le 2,
  $$
for some $s_0\in \mathbb{R}$, $0< \mathfrak{m}<\mathcal{M}<\infty$ and $0<\ell<1$, then analogously to Proposition~\ref{prop-area} there exists a positive constant $c_\ell$ depending only on $\ell$ such that for every $0\le k\le u\le2$ it holds
\begin{equation*}
\label{estimate-B-pme}c_\ell \mathfrak{m}(u-k)_+^2\le
\int_{0}^{(u-k)_{+}}\varphi^{\prime}(s+k)s\,ds\le \frac{\mathcal{M}}2 (u-k)^{\ell+1}_+.
\end{equation*}
Hence we will obtain a suitable energy inequality that will allow us to repeat the whole process and get H\"older continuity at each point. As in the singular case we can not obtain through this approach a better result for changing sign solutions, actually uniform H\"older regularity, since the constants in the energy inequality blow up near a degenerate point, where $\mathfrak{m}$ and $\mathcal{M}$ go to infinity for the rescaled problems, see Remark~\ref{rem-rem}.

\


\noindent{\large \textbf{Acknowledgments}}

\noindent All authors supported by the Spanish project MTM2014-53037-P.

\




\

\noindent\textbf{Addresses:}

\noindent\textsc{A. de Pablo: } Departamento de Matem\'{a}ticas, Universidad
Carlos III de Madrid, 28911 Legan\'{e}s, Spain. (e-mail: arturo.depablo@uc3m.es).

\noindent\textsc{F. Quir\'{o}s: } Departamento de Matem\'{a}ticas, Universidad
Aut\'{o}noma de Madrid, 28049 Madrid, Spain. (e-mail: fernando.quiros@uam.es).

\noindent\textsc{A. Rodr\'{\i}guez: } Departamento de Matem\'{a}tica Aplicada,  Universidad Polit\'{e}cnica de Madrid, 28040 Madrid, Spain.
(e-mail: ana.rodriguez@upm.es).

\end{document}